# Expansions of the solutions of the confluent Heun equation in terms of the incomplete Beta and the Appell generalized hypergeometric functions

C. Leroy[1] and A.M. Ishkhanyan[2]

[1]Laboratoire Interdisciplinaire Carnot de Bourgogne, CNRS UMR 6303, Université de Bourgogne, BP 47870, 21078 Dijon, France
[2]Institute for Physical Research of NAS of Armenia, 0203 Ashtarak, Armenia

*Corresponding author. Email: aishkhanyan@gmail.com

We construct several expansions of the solutions of the confluent Heun equation in terms of the incomplete Beta functions and the Appell generalized hypergeometric functions of two variables of the fist kind. The coefficients of different expansions obey four-, five-, or six-term recurrence relations that are reduced to ones involving less number of terms only in a few exceptional cases. The conditions for deriving finite-sum solutions via termination of the series are discussed.

**MSC numbers**: 33E30 Other functions coming from differential, difference and integral equations, 34B30 Special equations (Mathieu, Hill, Bessel, etc.), 30Bxx Series expansions

**Keywords:** Linear ordinary differential equation, confluent Heun equation, special functions, series expansions, recurrence relations

**1. Introduction**

The confluent Heun equation [1-3], the first confluent form of the general Heun equation [4], is widely encountered in mathematics, physics and engineering [1-2]. The special functions emerging from this equation generalize all the functions of the class of ordinary (i.e., not generalized) hypergeometric functions [5], as well as many other standard special functions, including the functions of the Bessel class, the Mathieu, the Whittaker-Hill, the spheroidal wave and Coulomb spheroidal wave functions [1-6].

The solutions of the Heun equations are written in terms of simpler mathematical functions only under rather restrictive conditions. Some of such solutions are derived by means of a rational change of the independent variable (generally accompanied with a transformation of the dependent variable also) [7-11], other reductions apply the Darboux transformation and other specific quadratures [12-13], however, the major known set is obtained by termination of infinite series solutions in terms of simpler special functions, in particular, functions of the hypergeometric class [14-35]. In the present paper we discuss the series solutions of the confluent Heun equation in terms of the incomplete Beta functions and the Appell generalized hypergeometric functions of two variables of the fist kind [5]. We show that in several cases these expansions also provide closed-form finite-sum solutions. It should be noted, however that in general the finite-sum solutions in terms of incomplete Beta functions are quasi-polynomials, so that the infinite series solutions in this case may suggest solutions of richer structure.

Expansions of the solutions of the confluent Heun equation in terms of the incomplete Beta functions have been first discussed in [25]. However, the expansions proposed in this reference apply to a rather restricted set of the involved parameters, namely, to the cases when the derivative of the solution of a confluent Heun equation is a solution of another confluent Heun equation with altered parameters [35]. In the present paper we follow a different approach. Applying an equation obeyed by a function involving the derivative of a solution of the confluent Heun equation, we show that incomplete Beta function expansions can be constructed for arbitrary parameters of the confluent Heun equation. As regards the expansions in terms of the Appell generalized hypergeometric functions of two variables of the first kind, such expansions have been previously suggested only for the general Heun equation [30,33]. For both expansion functions, the incomplete Beta and the Appell functions, we present two different types of expansions.

We discuss the conditions for deriving finite-sum solutions via termination of the series. The termination is possible only for one of the two types of the presented expansions. In this case, for the series governed by four-term recurrence relations for the expansion coefficients, these conditions result in a linear relation between the exponent parameters of the Heun equation and an additional polynomial equation for the accessory parameter. If the expansion coefficients are ruled by a five-term recurrence relation, one more restriction, which presents a polynomial equation for one of the exponent parameters, is imposed.

## 2. The first type of expansions

The confluent Heun equation is a second order linear differential equation having two regular singularities located at $z = 0$ and $1$, and an irregular singularity of rank 1 at $z = \infty$. We apply here the following form of this equation [3]:

$$\frac{d^2 u}{dz^2} + \left(\frac{\gamma}{z} + \frac{\delta}{z-1} + \varepsilon\right)\frac{du}{dz} + \frac{\alpha z - q}{z(z-1)} u = 0, \qquad (1)$$

which slightly differs from that adopted in [1-2] since in this form the parameters $\varepsilon$ and $\alpha$ are assumed to be independent. Note that in this notation the Whittaker-Ince limit [6] of the confluent Heun equation corresponds to the specification $\varepsilon = 0$. We denote the solution of Eq. (1) as $u = H_C(\gamma, \delta, \varepsilon; \alpha, q; z)$. It is assumed that this notation refers to a solution defined up to an arbitrary constant multiplier, not to one normalized to the unity at the origin.

We start with the observation that the function

$$v = z^{\gamma}(z-1)^{\delta} \frac{du}{dz} \quad (2)$$

obeys the following second order differential equation:

$$\frac{d^2 v}{dz^2} + \left( \frac{1-\gamma}{z} + \frac{1-\delta}{z-1} + \varepsilon - \frac{1}{z-z_0} \right) \frac{dv}{dz} + \frac{\Pi(z)}{z(z-1)(z-z_0)} v = 0, \quad (3)$$

where $z_0 = q/\alpha$ and $\Pi(z)$ denotes the quadratic polynomial

$$\Pi(z) = (\alpha - (\delta + \gamma - 1)\varepsilon)z^2 + (\gamma\varepsilon - z_0(2\alpha + \varepsilon(2 - \delta - \gamma)))z + z_0(z_0\alpha + \varepsilon - \gamma\varepsilon). \quad (4)$$

Equation (3) in general possesses three regular and an irregular singular points. Compared with Eq. (1), we have an additional regular singularity located at the point $z_0 = q/\alpha$. This is an apparent (non-logarithmic) singularity having characteristic exponents 0 and 2. The properties of such kind of regular singularities in connection with the Heun equations are recently studied by several authors [36-38]. Of course, the number of singularities reduces to three if this point coincides with one of the existing singularities of Eq. (1), that is if $z_0 = 0, 1, \infty$, which occurs if $q = 0$, $q = \alpha$ and $\alpha = 0$, respectively.

The expansions of the solutions of Eq. (1) in terms of the incomplete Beta functions and the Appell generalized hypergeometric functions are constructed in the following way. Consider a power-series expansion of a solution of Eq. (3) in the neighborhood of a point $z_1$ of the complex $z$-plane:

$$v = (z-z_1)^{\mu} \sum_{n=0}^{+\infty} a_n (z-z_1)^n . \quad (5)$$

The substitution of this series into Eq. (2) and subsequent integration results in the expansion

$$u = C_0 + \sum_{n=0}^{+\infty} a_n \left( \int z^{-\gamma}(z-1)^{-\delta}(z-z_1)^{\mu+n} dz \right). \quad (6)$$

In general the integrals involved in this expansion are expressed in terms of the Appell generalized hypergeometric functions of two variables of the first kind [5]. In several cases, however, the expansions are written in terms of simpler mathematical functions, in particular, in terms of the incomplete Beta functions. Note that the constant $C_0$ involved in this expansion is a significant component of the expansion; it is not an arbitrary constant. The value of this constant should be particularly specified in order to produce a valid expansion. This can be done by considering definite integration or by examining the limit $z \to 0$ when substituting the expansion into Eq. (1) (see examples below). If the integration is done over an interval between any two of singular points $0,1,2,\infty$, the integrals in Eq. (6) turn to Gauss hypergeometric functions $_2F_1$ [5].

The choice $z_1 = 0$ produces an immediate expansion in terms of the incomplete Beta functions ($|z| \leq 1$):

$$u = C_0 + \sum_{n=0}^{+\infty} a_n^{(0)} B(1 + n - \gamma + \mu, 1 - \delta; z), \quad \mu = 0, \gamma. \tag{7}$$

Substituting this expansion into Eq. (1) and taking the limit $z \to 0$ we readily find that here $C_0 = -\mu/q$ if $\operatorname{Re}(1 - \gamma + \mu) > 0$.

Similarly, choosing $z_1 = 1$ we get another expansion involving incomplete Beta functions with interchanged parameters as compared with the previous case:

$$u = C_0 + \sum_{n=0}^{+\infty} a_n^{(1)} (-1)^n B(1 - \gamma, 1 + n - \delta + \mu; z), \quad \mu = 0, \delta. \tag{8}$$

It is shown that in this case $C_0 = 0$ if $\operatorname{Re}(1 - \gamma) > 0$.

The coefficients of the constructed two expansions obey four-term recurrence relations. For instance, for the coefficients $a_n^{(0)}$ involved in Eq. (7) we have:

$$S_n a_n^{(0)} + R_{n-1} a_{n-1}^{(0)} + Q_{n-2} a_{n-2}^{(0)} + P_{n-3} a_{n-3}^{(0)} = 0, \tag{9}$$

where

$$S_n = q(n + \mu)(n - \gamma + \mu), \tag{10}$$

$$R_n = r_0 + r_1(n + \mu) + r_2(n + \mu)^2, \quad Q_n = q_0 + q_1(n + \mu) + q_2(n + \mu)^2, \tag{11}$$

$$P_n = \alpha(\alpha + \varepsilon(1 + n - \gamma - \delta + \mu)), \tag{12}$$

where $\mu = 0, \gamma$ and the parameters $r_{0,1,2}$, $q_{0,1,2}$ do not depend on $n$. It is seen that the recurrence relation becomes three-term if $q = 0$ ($S_n \equiv 0$) or if $\alpha = 0$ ($P_n \equiv 0$). This is of course an expected result since in these cases Eq. (3) presents another confluent Heun equation with altered parameters [35].

If $q \neq 0$, the series is left-hand side terminated at $n = 0$ if $S_0 = 0$, i.e., if $\mu = 0$ or $\mu = \gamma$. It will terminate from the right-hand side if three successive coefficients vanish for some $N = 1, 2, \ldots$: $a_N \neq 0$, $a_{N+1} = a_{N+2} = a_{N+3} = 0$. From the equation $a_{N+3} = 0$ we find that the termination is possible if $P_N = 0$. For non-zero $\alpha$ this is the case if

$$\alpha = -\varepsilon(1 + N - \gamma - \delta + \mu), \quad \mu = 0, \gamma. \tag{13}$$

Since $\gamma, \delta$ determine the characteristic exponents of the regular singularities $z = 0,1$ and the ratio $\alpha/\varepsilon$ characterizes the irregular singularity $z = \infty$ [2], this equation can be viewed as a linear relation between the exponent parameters of the Heun equation). Further, it is verified by direct inspection that if this relation holds, the remaining two equations, $a_{N+1} = 0$ and $a_{N+2} = 0$, become linearly dependent so that only one additional restriction is imposed on the parameters of the confluent Heun equation. It is readily shown that this restriction presents a polynomial equation of the degree $N+1$ for the accessory parameter $q$.

As already stated above, in the general case $z_1 \neq 0,1$ Eq. (6) presents an expansion in terms of the Appell generalized hypergeometric functions of the first kind:

$$u = C_0 + \sum_{n=0}^{+\infty} a_n^{(z_1)} \frac{(-1)^\delta}{(-z_1)^{-\mu-n}} \frac{z^{1-\gamma}}{1-\gamma} F_1\left(1-\gamma;\ \delta, -\mu-n;\ 2-\gamma;\ z, \frac{z}{z_1}\right). \tag{14}$$

From Eq. (3) we get that $\mu = 0$ or $2$ if $z_1 = z_0 = q/\alpha$. We note that it follows from the properties of the Appell functions that in all above cases the values adopted by the expansion functions $u_n$ at the point $z = 1$ can be written in terms of the hypergeometric functions [5]:

$$u_n(z=1) \sim F_1(a_1; b_1, b_2; a_1+1; z, y)\big|_{z \to 1} = \frac{\Gamma(a_1+1)\Gamma(1-b_1)}{\Gamma(a_1+1-b_1)} {}_2F_1(a_1, b_2; a_1+1-b_1; y).$$

(By a rational transformation of $z$, a similar result is achieved for the points $z = z_0, \infty$ as well). This is of course well seen if one applies definite integration in deriving Eq. (6).

In several cases the Appell functions involved in this expansion are also reduced to incomplete Beta functions. For instance, this occurs if $\delta = 0$:

$$u(\delta = 0) = C_0 + \sum_{n=0}^{+\infty} a_n^{(z_1)}(-z_1)^n B\left(1-\gamma,\ 1+n+\mu;\ \frac{z}{z_1}\right), \tag{15}$$

or if $\gamma = 0$:

$$u(\gamma = 0) = C_0 + \sum_{n=0}^{+\infty} a_n^{(z_1)}(1-z_1)^n B\left(1+n+\mu,\ 1-\delta;\ \frac{z-z_1}{1-z_1}\right). \tag{16}$$

In general the coefficients $a_n^{(z_1)}$ of Eq. (14) obey a five-term recurrence relation:

$$T_n a_n^{(z_1)} + S_{n-1} a_{n-1}^{(z_1)} + R_{n-2} a_{n-2}^{(z_1)} + Q_{n-3} a_{n-3}^{(z_1)} + P_{n-4} a_{n-4}^{(z_1)} = 0, \tag{17}$$

where

$$T_n = z_1(z_1-1)(z_1-z_0)(n+\mu)(n-1+\mu), \tag{18}$$

$$S_n = s_1 \mu_1 + s_2 \mu_1^2,\ \ R_n = r_0 + r_1 \mu_1 + r_2 \mu_1^2,\ \ Q_n = q_0 + q_1 \mu_1 + \mu_1^2,\ \ \mu_1 = n+\mu, \tag{19}$$

$$P_n = \alpha + \varepsilon(1+n-\gamma-\delta+\mu). \tag{20}$$

Here again the parameters $s_{1,2}$, $r_{0,1,2}$, $q_{0,1}$ do not depend on $n$. This recurrence relation involves just four successive terms only if $z_1 = 0, 1, \infty, z_0$. The recurrence relation for $z_1 = 0$ is given by Eq. (9). For $z_1 = 1$ we get a similar recurrence relation. In the case when $z_1 = z_0$ we have $T_n = 0$ and $S_n = z_0(z_0 - 1)(n + \mu)(n - 2 + \mu)$, so that in this case the series may left-hand side terminate if $\mu = 0$ or $\mu = 2$. If $z_1 \neq 0, 1, z_0$ the series may left-hand side terminate if $\mu = 0$ or $\mu = 1$. As stated above, the singularity $z = z_0$ is an apparent one, so that in all cases consistent power-series expansions, without logarithmic terms, are constructed for both exponents, not only for the greater one.

For non-zero $\alpha$ the series (14) may terminate from the right-hand side at some $N = 1, 2, \ldots$ if $P_N = 0$, that is if

$$\alpha = -\varepsilon(1 + N - \delta - \gamma + \mu), \tag{21}$$

and additionally if $a_{N+1} = a_{N+2} = a_{N+3} = 0$ (in the general case of the five-term relation) or $a_{N+1} = a_{N+2} = 0$ (in the case of a four-term relation). Though in general one could expect that these equations impose respectively three or two additional restrictions on the parameters of the confluent Heun equation, however, the close inspection of coefficients (18)-(20) reveals that if Eq. (21) holds these equations become dependent, so that less restrictions are imposed: just two equations in the case of the five-term relation and only one if the relation involves four-terms. In the four-term case this is a polynomial equation of the degree $N + 1$ for the accessory parameter $q$. Correspondingly, for given $N$, in this case the termination occurs for $N + 1$ values of $q$. In the five-term case we have two polynomial equations each of the degree $N$ if $\mu = 0$ and of the degrees $N + 1$ and $N + 2$ if $\mu = 2$. Correspondingly, for given $N$, in the case of the five-term recurrence relation the termination is achieved for $N^2$ or $(N + 1)(N + 2)$ particular sets of pairs $\{q, \varepsilon\}$ for $\mu = 0$ or $\mu = 2$, respectively.

### 3. The second type of expansions

It is possible to construct other expansions in terms of the incomplete Beta functions as well as in terms of the Appell generalized hypergeometric functions. For instance, this can be done if a preliminary change of the dependent variable is applied. To demonstrate the approach, suppose $\varepsilon \neq 0$ and consider the transformation $u = e^{-\varepsilon z/2} w(z)$, which reduces Eq. (1) to the equation

$$\frac{d^2w}{dz^2} + \left(\frac{\gamma}{z} + \frac{\delta}{z-1}\right)\frac{dw}{dz} + \frac{\Pi(z)}{z(z-1)}w = 0, \tag{22}$$

where $\Pi(z)$ is the following polynomial:

$$\Pi(z) = \frac{1}{4}\left(-\varepsilon^2 z^2 + (4\alpha - 2(\gamma+\delta)\varepsilon + \varepsilon^2)z + 2\gamma\varepsilon - 4q\right). \tag{23}$$

For non-zero $\varepsilon$ this is a quadratic polynomial that can be factorized as $p_0(z-z_1)(z-z_2)$ with $p_0 = -\varepsilon^2/4$. It is then readily verified that the function

$$v = z^\gamma (z-1)^\delta \frac{dw}{dz} \tag{24}$$

obeys the following second order equation:

$$\frac{d^2v}{dz^2} + \left(\frac{1-\gamma}{z} + \frac{1-\delta}{z-1} - \frac{1}{z-z_1} - \frac{1}{z-z_2}\right)\frac{dv}{dz} + \frac{p_0(z-z_1)(z-z_2)}{z(z-1)}v = 0. \tag{25}$$

In general, this equation has five singular points; we have two additional singularities located at $z = z_1$ and $z = z_2$. As in the case of Eq. (3), these are apparent singularities having characteristic exponents 0, 2. Of course, the number of singularities decreases if one or both of these singularities coincide with already existing real singularities $z = 0, 1$.

Now again considering a power-series expansion of a solution of Eq. (25) in the neighborhood of a point $z_i$ of the complex $z$-plane:

$$v = (z-z_i)^\mu \sum_{n=0}^{+\infty} a_n (z-z_i)^n, \tag{26}$$

we arrive at the expansion (compare with Eq. (6))

$$u = e^{-\varepsilon z/2}\left[C_0 + \sum_{n=0}^{+\infty} a_n \left(\int z^{-\gamma}(z-1)^{-\delta}(z-z_i)^{\mu+n}dz\right)\right]. \tag{27}$$

The integrals here are of the same form as the ones involved in Eq. (6), hence, similar developments as the above expansions Eqs. (7), (8) and Eq. (14) apply. Indeed, the choice $z_i = 0$ produces the expansion

$$u = e^{-\varepsilon z/2}\left[C_0 + \sum_{n=0}^{+\infty} a_n^{(0)} B(1+n-\gamma+\mu, 1-\delta; z)\right], \quad \mu = 0, \gamma, \tag{28}$$

where $C_0 = -\mu/(q-\gamma\varepsilon/2)$ if $\text{Re}(1-\gamma+\mu) > 0$, and the choice $z_i = 1$ produces another incomplete Beta function expansion written as

$$u = e^{-\varepsilon z/2}\left[C_0 + (-1)^\delta \sum_{n=0}^{+\infty} a_n^{(1)}(-1)^{-n-\mu} B(1-\gamma, 1+n-\delta+\mu; z)\right], \quad \mu = 0, \delta, \tag{29}$$

where $C_0 = 0$ if $\text{Re}(1-\gamma) > 0$.

Finally, in the general case $z_i \neq 0,1$ we have a different expansion in terms of the Appell generalized hypergeometric functions of the first kind:

$$u = e^{-\varepsilon z/2} \left[ C_0 + \sum_{n=0}^{+\infty} a_n^{(z_1)} \frac{(-1)^{\delta}}{(-z_i)^{-\mu-n}} \frac{z^{1-\gamma}}{1-\gamma} F_1\left(1-\gamma; \delta, -\mu-n; 2-\gamma; z, \frac{z}{z_i}\right) \right]. \quad (30)$$

Though the above two types of expansions, Eqs. (7)-(8),(14) and (28)-(30), have much in common, there are significant differences between the expansions of the first and the second types that should be mentioned. The differences concern the recurrence relations between the successive coefficients of the expansions. First, for the first type expansions the recurrence relation in general is five-term, while for the second type expansions the relation in general is six-term that may involve fewer terms only under rather restrictive conditions. Furthermore, the presented second type expansions do not terminate from the right-hand side.

Indeed, the coefficients $a_n^{(0)}$ of expansion (28) obey the following six-term relation:

$$K_n a_n^{(0)} + T_{n-1} a_{n-1}^{(0)} + S_{n-2} a_{n-2}^{(0)} + R_{n-3} a_{n-3}^{(0)} + Q_{n-4} a_{n-4}^{(0)} + P_{n-5} a_{n-5}^{(0)} = 0 \quad (31)$$

where the first and the last coefficients read

$$K_n = -z_1 z_2 (\mu+n)(\mu+n-\gamma), \quad P_n = -\frac{\varepsilon^2}{4}. \quad (32)$$

It is seen that this recurrence relation may reduce to one involving five or less successive terms only if $z_1 z_2 = 0$. Furthermore, since it is supposed that $\varepsilon \neq 0$, from the second equation (32) we conclude that the series cannot terminate from the right-hand side. If non of $z_1, z_2$ is zero, the series is left-hand side terminated at $n = 0$ if $K_0 = 0$, i.e., if $\mu = 0$ or $\mu = \gamma$.

## 4. Summary

Thus, we have presented several expansions of the solutions of the confluent Heun equation in terms of the incomplete Beta functions and the Appell generalized hypergeometric functions of two variables of the first kind. We have seen that the coefficients of the expansions obey four-, five-, or six-term recurrence relations. Discussing the conditions for the series to terminate, we have seen that this is possible only for the first type expansions for which the recurrence relations for the expansion coefficients involve four or five terms. In the cases when the series are governed by four-term relations, the conditions for termination result in a relation between the exponent parameters of the confluent Heun equation and a polynomial equation for the accessory parameter. If the expansion coefficients are ruled by a five-term recurrence relation, one more polynomial equation should be satisfied by the involved parameters.

It is understood that other expansions can be constructed if a preliminary change of the dependent and/or independent variables is applied. It is a question if such a transformation is able to produce simpler recurrence relations. A final remark is that other expansions of the discussed type as well as expansions in terms of other higher transcendental functions, e.g., the Goursat generalized hypergeometric functions [35] can be suggested if equations for products of the Heun functions [11,13,39], or Laplace transform and its inverse [40,41], or an integral Euler transformation [42], or other integral equations, relations, and representations [1-3,43-49] are applied.

## Acknowledgments

This research has been conducted within the scope of the International Associated Laboratory (CNRS-France & SCS-Armenia) IRMAS. The research has received funding from the European Union Seventh Framework Programme (FP7/2007-2013) under grant agreement No. 295025 – IPERA. The work has been supported by the Armenian State Committee of Science (SCS Grant No. 13RB-052).